\documentclass[12pt]{article}
\usepackage[cmtip,arrow]{xy}
\usepackage{pb-diagram,pb-xy}
\usepackage{amssymb} \usepackage{amsfonts}
\usepackage{amsmath}
\usepackage{amscd}
\usepackage{url}
\usepackage{amsthm}

\newcommand{\F}{{\mathbb F}}
\newcommand{\Z}{{\mathbb Z}}

\newcommand{\R}{{\mathbb R}}
\newcommand{\C}{{\mathbb C}}

\newcommand{\N}{{\mathbb N}}
\begin{document}
\title{Examples of non connective $C^{\star}$-algebras}
\author{A. G\c{a}sior\thanks{The first author is supported by the Polish National Science Center grant DEC2017/01/X/ST1/00062.}, A. Szczepa\'nski}
\date{\today}
\maketitle
\begin{abstract}
We give an example of two infinite families of not connective groups. Both of them are generalization
of the 3-dimensional Hantzsche - Wendt group.
\end{abstract}
\vskip 10mm

\section{Introduction}

For a Hilbert space ${\cal H},$ we denote by $L({\cal H})$
the $C^{\star}$-algebra of bounded and linear operators on ${\cal H}.$
The ideal of compact operators is denoted by ${\cal K}\subset L({\cal H}).$
For the $C^{\star}$ - algebra $A,$ the cone over $A$ is defined as $CA =C_{0}[0, 1)\otimes A,$ the suspension of $A$ as
$SA =C_{0}(0, 1)\otimes A.$

\newtheorem{defi}{Definition}
\begin{defi}
Let $A$ be a $C^{\star}$-algebra and $n\in\N, n\geq 1.$ $A$ is {\em connective} if there is a $\star$-monomorphism
$$\Phi:A\to\prod_{n} CL({\cal H})/\bigoplus_{n} CL({\cal H})$$
which is liftable to a completely positive and contractive map $\phi:A\to \prod_{n} CL({\cal H}).$
\end{defi}
\vskip 1mm
\noindent
For a discrete group $G,$ we define $I(G)$ to be the augmention ideal, i.e. the kernel of the trivial representation
$C^{\star}(G)\to\C.$
$G$ is called connective if $I(G)$ is a connective $C^{\star}$-algebra.
From definition (see \cite[p. 4921]{Dad}) connectivity
of $G$ may be viewed as a stringent topological property that accounts simultaneously for the quasidiagonality of
$C^{\star}(G)$ and the verification of the Kadison-Kaplansky conjecture for certain classes of groups.
Here we can referring to conjecture from 2014 \cite[p. 166]{D}. {\em If $G$ is a discrete, countable, torsion-free, amenable
group, then the natural map
$$[[I(G),\mathcal{K}]]\to KK(I(G),\mathcal{K})\cong K^{0}(I(G))$$ is an isomorphism of groups.}
Where $KK(I(G),\mathcal{K})$ is the Kasparov group and $[[I(G),\mathcal{K}]]$ is a group of the homotopy classes of
asymptotic morphisms. In 2017 M. Dadarlat found an amenable and not connective group $G_2$ for which the above conjecture fails \cite[Cor. 3.2]{D2}.

Connective groups must be torsion-free, \cite[Remark 2.8 and 4.4]{D3}.
Here is a short list of such groups:
\begin{itemize}
\item[1.] a countable torsion free nilpotent groups, \cite[Th.4.3]{D3};
\item[2.] let $0\to N\to G\to H\to 0$ be a central extension of
discrete countable amenable groups where $N$ is torsion-free. If $H$ is connected then so does $G$;
        \cite[Th. 4.1]{D3};
\item[3.] wreath product of connected groups is a connected group \cite[Th.3.2]{DPS};
\item[4.] a torsion-free crystallographic group is connective if and only if is locally indicable if and only if is diffuse (see below) and \cite{D2}.
\end{itemize}
A discrete group $G$ is called locally indicable if every finitely generated non-trivial subgroup $L$ of $G$ has an infinite abelianization.
The group $G$ is called diffuse if every non-empty finite subset $A$ of $G$ has an element $a\in A$ such that for any $g\in G,$ either
$ga$ or $g^{-1}a$ is not in $A,$ \cite{D2}, \cite{GLS}.
\vskip 1mm
\noindent
More examples of nonabelian connective groups were exhibited in \cite{D2}, \cite{Dad}, \cite{DPS}.
\vskip 1mm
The above group $G_2$ is a 3-dimensional, torsion-free crystallographic group, where
a crystallographic group $\Gamma,$ of dimension $n$ is a cocompact and discrete subgroup of the isometry group $E(n) = O(n)\ltimes\R^n$
of the Euclidean space $\R^n$. $\Gamma$ is cocompact if and only if the orbit space $E(n)/\Gamma$ is compact.
\noindent
From Bieberbach theorems (see \cite[Chapter 1]{S}) any crystallographic group $\Gamma$ defines a short exact sequence
$$0\to\Z^n\to \Gamma\to H\to 0,$$
where a free abelian group $\Z^n$ is a maximal abelian subgroup and $H$ is a finite group.
$H$ is sometimes called a holonomy group of $\Gamma.$ The above group $G_2$ is isomorphic to the subgroup $E(3)$
and is generated by
{\scriptsize $$G_{2}\cong \text{gen}\{A =(\left[
    \begin{array}{clcr}
    1 & 0 & 0\\
    0 & -1 & 0\\
    0 &  0 & -1
    \end{array}
    \right],(1/2,1/2,0)),
    B = (\left[
    \begin{array}{clcr}
    -1 & 0 & 0\\
    0 & 1 & 0\\
    0 & 0 & -1
    \end{array}
    \right],(0,1/2,1/2))\}.
    $$}
\newtheorem{prop}{Proposition}
\vskip 2mm
\noindent
A torsion-free crystallographic group is called a Bieberbach group. The orbit space $\R^n/\Gamma$ of a Bieberbach
group is a $n$-dimensional closed flat Riemannian manifold $M$ with holonomy group isomorphic to $H.$
\vskip 1mm
\noindent
A general characterization of connective Bieberbach groups is given in \cite{D2}.
The following two theorems give us a landscape of them.
\newtheorem{theo}{Theorem}
\begin{theo}\label{main} {\em (\cite[Theorem 1.2]{D2})}
Let $\Gamma$ be a Bieberbach group. The following assertions are equivalent.
\begin{itemize}
\item[1.] $\Gamma$ is connective
\item[2.] Every nontrivial subgroup of $\Gamma$ has a nontrivial center.
\item[3.] $\Gamma$ is a poly-$\Z$ group
\item[4.] $\widehat{G}\setminus\{\iota\}$ has no nonempty compact open subsets.
\end{itemize}
\noindent
The unitary dual $\widehat{G}$ of $G$ consists of equivalence classes of irreducible unitary representations of $G.$
$\iota$ denotes the trivial representation.
\end{theo}
\begin{theo}\label{main1} {\em (\cite[Theorem 1.1]{D2})}
A Bieberbach group with a finite abelianization is not connective.
\end{theo}

In our note we give an example of two infinite families of not connective groups. Both of them are generalization
of the 3-dimensional Hantzsche-Wendt group $G_2.$

\section{Examples}

\newtheorem{exe}{Example}
\begin{exe}\label{ex1}
{\em (\cite[Definition 9.1]{S})}Let $n\geq 3$. By generalized Hantzsche-Wendt (GHW for short) group we shall understand any torison-free crystallographic groups of rank $n$ with a holonomy group $(\Z_2)^{n-1}.$
\end{exe}
\vskip 2mm
\begin{exe}\label{ex2}
{\em (\cite[Definition]{CL}, \cite[Definition 1]{S1})} Let $n\geq 0.$
A group
$$G_n = \{x_1,x_2,...,x_n\mid x_{i}^{-1}x_{j}^{2}x_{i}x_{j}^{2}, \forall i\neq j\}$$
we shall call a combinatorial Hantzsche-Wendt group.
\end{exe}
For the properties of GHW groups we refer to \cite[Chapter 9]{S}.
We have $G_0 = 1$ and $G_1\simeq\Z.$
Combinatorial Hantzsche-Wendt groups are torsion-free, see \cite[Theorem 3.3]{CL} and for $n\geq 2$ are nonunique
product groups. A group $G$ is called a unique product group if given two nonempty finite subset $X,Y$ of $G,$
there exists at least one element $g\in G$ which has a unique representation $g=xy$ with $x\in X$ and $y\in Y.$
We are ready to present our main result.
\begin{prop}
Generalized Hantzsche-Wendt groups with trivial center and nonabelian, combinatorial Hantzsche-Wendt groups are not connective.
\end{prop}
\noindent
{\bf Proof:} From \cite[Remark 2.8 (i)]{D3} the connectivity property is inherited by subgroups. Let $G$ be any group from family of GHW groups
or family of combinatorial Hanzsche-Wendt groups. In both cases a group $G_2$ is a subgroup of $G$. In the first case
it follows from \cite[Proposition 9.7]{S}. In the second case it follows from definition, see \cite[Prop. 3.4]{CL}.
\vskip 1mm
\noindent
Note that in the case of GHW groups we can also use Theorem \ref{main1},
since all these groups have a finite abelianizations.
\vskip 5mm \hskip 120mm $\Box$
\vskip 1mm
\noindent

\newtheorem{rem}{Remark}
\begin{rem}
From {\em \cite{S1}}, for $n\geq 3, G_n$ has a non-abelian free subgroup. Hence is not amenable.
\end{rem}
\begin{rem}
The counterexample to the Kaplansky unit conjecture was given in 2021 by {\em G. Gardam \cite{Ga}}.
It was found in the group ring $\F_2[G_2].$
The Kaplansky unit conjecture states that every unit in $K[G]$ is of the form $kg$ for
$k\in K \setminus \{0\}$ and $g \in G.$
\end{rem}
\vskip 5mm
\noindent
{\bf Acknowledgements} We thank the referee for a number of suggestions that improved the exposition.

\vskip 2mm
\noindent
Institute of Mathematics, Maria Curie-Sk{\l}odowska University\\
Pl. Marii Curie-Sk{\l}odowskiej 1\\
20-031 Lublin\\
Poland\\
E-mail: \texttt{anna.gasior@poczta.umcs.lublin.pl}
\vskip 5mm
\noindent
Institute of Mathematics, University of Gda\'nsk\\
ul. Wita Stwosza 57,\\
80-952 Gda\'nsk,\\
Poland\\
E-mail: \texttt{matas@univ.gda.pl}
\end{document}